\documentclass[final,5p,times]{elsarticle}

\usepackage{amsfonts}
\usepackage{amsmath}
\usepackage{amssymb}
\usepackage{graphicx}

\def\ooo{{\Omega}}
\def\<{\left<}
\def\>{\right>}
\def\({\left(}
\def\){\right)}
\def\ff{\forall }
\def\D{\Delta }
\def\9{{\infty}}
\def\barr{\begin{array}}
\def\earr{\end{array}}
\def\ov{\overline}
\def\wt{\widetilde}
\def\vp{{\varepsilon}}

\def\pp{{\partial}}
\def\dd{\displaystyle}
\def\vf{{\varphi}}
\def\lbb{{\lambda}}
\def\g{{\gamma}}
\def\a{{\alpha}}
\def\b{{\beta}}
\def\pas{\mathbb{P}\mbox{-a.s.}}
\def\de{{\delta}}
\def\cala{{\mathcal{A}}}
\def\calf{{\mathcal{F}}}
\def\calo{\mathcal{O}}
\def\calf{\mathcal{F}}
\def\calh{\mathcal{H}}
\def\calx{\mathcal{X}}
\def\vsp{\vspace*{1,5mm}\\ }
\def\cala{{\mathcal{A}}}
\def\vspp{\vspace*{1,5mm}\\ }
\def\n{\noindent }

\def\inr{\mbox{ in }}
\def\onr{\mbox{ on }}
\def\D{{\Delta}}
\def\rr{{\mathbb{R}}}
\def\dd{\displaystyle}
\def\Res{{\rm Re}}
\def\Ims{{\rm Im}}
\def\S{Stratonovich}
\def\3{\subset }
\def\na{{\nabla}}

\newtheorem{theorem}{Theorem}
\newtheorem{remark}[theorem]{Remark}

\journal{Mathematics....}

\begin{document}

\begin{frontmatter}
\title{Internal stabilization of the Oseen--Stokes equations by Stratonovich noise}
\author{Viorel Barbu}
\ead{vb41@uaic.ro}
\address{"Al.I. Cuza" University and "Octav Mayer" Institute of Mathematics of Romanian Academy, Ia\c si, Romania}
\begin{abstract}
One designs a Stratonovich noise feedback controller with support in an arbitrary open subset $\calo_0$ of $\calo$ which exponentially stabilizes in probability, that is with probability one, the Oseen--Stokes systems in a domain $\calo\3\rr^d$, $d=2,3$. This completes the stabilization results from the author's work \cite{5} which is concerned with design of an Ito noise stabilizing controller.
\end{abstract}

\end{frontmatter}

\section{Introduction and statement of the problem}

This work is concerned with internal stabilization via Stra\-to\-no\-vich noise feedback controller of Oseen--Stokes system
\begin{equation}\label{e1}
\barr{l} \dd\frac{\pp X}{\pp t}-\nu\D
X+(f\cdot\na)X+(X\cdot\na)g=\na p\ \inr\ (0,\9)\times\calo,\vspp
\na\cdot X=0\quad\inr\ (0,\9)\times\calo,\vspp X(0,\xi)=x(\xi),\
\xi\in\calo,\quad X=0\onr (0,\9)\times\pp\calo.
\earr\end{equation} Here, $\calo$ is an open and bounded  subset
of $\rr^d$, $d=2,3$, with smooth boundary $\pp\calo$ and $f,g\in
C^2(\ov\calo;\rr^d)$ are given functions. In the special case
$g\equiv0$, system \eqref{e1} describes the dynamic of a fluid
Stokes flow with partial inclusion of convection
acce\-le\-ra\-tion $(f\cdot\na)X$ $(X$ is the velocity field). The
same equation describes the disturbance flow induced by a moving
body with velocity $f$ through the fluid. Should we mention also
that in the special case $f\equiv g\equiv X_e$, where $X_e$ is the
equilibrium (steady-state) solution of the Navier--Stokes equation
\begin{equation}\label{e2}
\barr{l} \dd\frac{\pp X}{\pp t}-\nu\D X+(X\cdot\na)X=\na
p+f_e,\vspp \na\cdot X=0,\quad X|_{\pp\calo}=0,\earr\end{equation}
and $f_e\in C(\ov\calo;\rr^d)$, system \eqref{e1} is the
linearization of \eqref{e2} around $X_e$. In this way, the
stabilization of \eqref{e1} can be interpreted as the first order
stabilization procedure of steady-state Navier--Stokes
flows.\smallskip

Our aim here is to design a stochastic feedback controller of the form
\begin{equation}\label{e3}
u={\bf1}_{\calo_0}\sum^M_{k=1}R_k(X)\circ\dot\b_k,\
  R_k\in L((L^2(\calo))^d),\end{equation}
which stabilizes in probability system \eqref{e1} and has the
support in an open subdomain $\calo_0$ of $\calo$.

Here, $\{\b_k\}^M_{k=1}$ is a system of mutually independent
Brownian motions on a pro\-ba\-bi\-lity space
$\{\ooo,\calf,\mathbb{P}\}$ with filtration $\{\calf_t\}_{t>0}$
while the cor\-res\-pon\-ding closed loop system
\begin{equation}\label{e4}
\barr{l} dX-\nu\D X\,dt+(f\cdot\na)X\,dt+(X\cdot\na)g\,dt\vsp
\qquad\hfill={\bf1}_{\calo_0}\dd\sum^M_{k=1}R_k(X)\circ d\b_k+\na
p\,dt\vsp X(0)=x\earr\end{equation} is taken in \S\ sense (see,
e.g., \cite{1})  and this is the significance of the symbol
$R_k(X)\circ\dot\b_k$ in the expression of the noise controller
\eqref{e3}. In the following, the terminology Stratonovich feedback controller or Ito feedback controller refer to the sense in which the corresponding stochastic equation \eqref{e4} is con\-si\-dered: in the Stratonovich sense or, respectively, the Ito sense. We have denoted by ${\bf1}_{\calo_0}$ the
characteristic function of the  open set
$\calo_0\subset\calo$.\smallskip

In \cite{4}, \cite{5}, \cite{6}, \cite{7},  the author has
designed similar sta\-bi\-li\-zable Ito noise controllers for
equation \eqref{e1} and related Navier--Stokes equations. However,
it should be said that, with respect to Ito noise controllers, the
\S\ feedback   controller \eqref{e3} has the advantage to be
stable with respect to smooth changes $\dot\b_k^\vp$ of the noise
$\dot\b_k$ and this fact is crucial not only from the conceptual
point of view, but also for numerical simulations and practical
implementation into   system \eqref{e1} of the random stabilizable
feedback controller
$$u_\vp(t)={\bf1}_{\calo_0}\sum^M_{k=1}R_k(X(t))\dot\b^\vp_k(t),$$
where $\dot\b^\vp_k$ is a smooth approximation of
$\b^\vp_k$. If instead \eqref{e3} we take $u$ to be an Ito stabilizable feedback controller, then the corresponding Ito stochastic closed loop equation is convergent for $\vp\to0$ to a Stratonovich equation of the form \eqref{e4} which might be unstable because the stabilization effect of the noise controller is given by the Ito's formula which is valid for Ito stochastic equation only. Thus, for numerical implementations of a stabilizanble noise feedback controller $u$ it is essential that it is of Stratonovich type. It should be said, however, that the option for an Ito stabilizable noise controller, like in \cite{4}--\cite{7}, or a Stratonovich one, as in this work, is function of the specific techniques used to insert such a noise controller into system \eqref{e1}; by direct simulation or by numerical approximating procedure.\smallskip

As regards the literature on  stabilization of linear differential
systems by \S\ noise, the pioneering works \cite{2}, \cite{3}
should be primarily cited. For linear PDEs, this procedure was
developed in \cite{8}, \cite{9} which are related to this work.
For ge\-ne\-ral results on internal stabilization of
Navier--Stokes systems with deterministic feedback controllers, we
refer to \cite{4a}. (See also \cite{11} for noise stabilization effect into a different nonlinear PDE context.)

\section{The noise stabilizing feedback controller}

Consider the standard space of free divergence vectors\break $H=\{y\in(L^2(\calo))^d;$ $\na\cdot y=0$ in $\calo$, $y\cdot n=0$ on $\pp\calo\}$ and denote by $\cala_0:D(\cala_0)\3H\to H$ the realization of the Oseen--Stokes operator in this space, that is,
\begin{equation}\label{e5}
\cala_0y=P(-\nu\D y+(f\cdot\na)y+(y\cdot\na)g),\quad y\in
D(\cala_0),\end{equation} where
$D(\cala_0)=H\cap(H^2(\calo))^d\cap(H^1_0(\calo))^d.$ Here, $P$ is
the Leray projector on $H$ and $H^2(\calo),$ $H^1_0(\calo)$ are
standard Sobolev spaces on $\calo$. In the following, it will be
more convenient to represent equation \eqref{e1} in the complex
Hilbert space $\calh=H+iH$ by extending $\cala_0$ to
$\cala:D(\cala)\3\calh\to\calh$ via standard procedure,
$\cala(x+iy)=\cala_0x+i\cala_0y$. The operator $\cala$ has a
countable set of eigenvalues $\{\lbb_j\}^\9_{j=1}$ (eventually
complex) with the cor\-res\-pon\-ding eigenvectors $\vf_j$. Denote by
$\cala^*$ the adjoint of $\cala$ with eigenvalues $\ov\lbb_j$ and
eigenvectors $\vf^*_j$. Each eigenvalue is repeated in the
following according to its algebraic multiplicity $m_j$. Normalizing the system $\{\vf_j\}^\9_{j=1}$, we see that
$$|\nabla\vf_j|^2_\calh=\lbb_j-\left<(f\cdot\nabla)\vf_j+(\vf_j\cdot\nabla)g,\vf_j\right>,\ \ff j,$$and since, by the Fredholm--Riesz theory, $|\lbb_j|\to+\9$ as $j\to\9$, we infer that ${\rm Re}\ \lbb_j\to+\9$ as $j\to\9$. We denote
by $N$ the minimal number of eigenvalues $\lbb_j$ for which
\begin{equation}\label{e6}
\Res\,\lbb_j>0\mbox{ for }j\ge N,\
\lbb_1+\lbb_2+\cdots+\lbb_N>0.\end{equation} (In the above
sequence, each $\lbb_j$ is taken together its conjugate
$\ov\lbb_j$ and, clearly, there is such a natural number $N$.)

Set $\calx_u{=}{\rm lin\, span}\{\vf_j\}^N_{j=1}$ and denote by
$\calx_s$ the algebraic complement of $\calx_u$ in $\calx$. It is
well known that $\calx_u$ and $\calx_s$ are both invariant for
$\cala$ and, if we set
$$\cala_u=\cala|_{\calx_u},\quad\cala_s=\cala|_{\calx_s},$$
we have for their spectra $\sigma(\cala_u)=\{\lbb_j\}^N_{j=1},$
$\sigma(\cala_s)=\{\lbb_j\}^\9_{j=N+1}$ and, since $-\cala_s$ is
the generator of an analytic $C$-semigroup $e^{-\cala_st}$ in
$\calx_s$, we have
\begin{equation}\label{e7}
\|e^{-\cala_st}\|_{L(\calh)}\le C\exp(-\Res\,\lbb_{N+1}t),\quad
t\ge0,\end{equation} (see, e.g., \cite{4}, p. 14). In the
following, we shall assume that
\begin{itemize}
\item[(i)] {\it All the eigenvalues $\{\lbb_j\}^N_{j=1}$ are semisimple.}
\end{itemize}
This means that  the algebraic multiplicity of each $\lbb_j$,
$j=1,...,N$, coincides with its geometric multiplicity or, in
other words, the finite-dimensional operator (matrix) $\cala_u$ is
diagonalizable. As we will see later on, this assumption is not
essentially necessary but it simplifies however the construction
of the stabilizing controller because it reduces the unstable part
of the system to a diagonal finite-dimensional differential
system. In particular, it follows by (i) that we can choose the
dual systems $\{\vf_j\}$ and $\{\vf^*_j\}$ in such a way that
\begin{equation}\label{e8}
\<\vf_i,\vf^*_j\>=\de_{ij},\quad i,j=1,...,N.\end{equation}
(Here, and everywhere in the following, $\<\cdot,\cdot\>$ stands for the scalar product in $\calh$ and $H$. By $|\cdot|_\calh$ and $|\cdot|_H$ we denote the corresponding norms.)

We note that the uncontrolled Oseen--Stokes system \eqref{e1} can be rewritten in the space $\calh$ as
\begin{equation}\label{e9}
\frac{dX}{dt}+\cala X=0,\  t\ge0,\quad X(0)=x,\end{equation}
and setting $X_u=\sum^N_{j=1}y_j\vf_j,$ $X_s=(I-P_N)X$, where $P_N$ is the algebraic projector on $\calx_u$, we have
\begin{eqnarray}
\dd\frac{dX_u}{dt}+\cala_uX_u=0,\quad X_u(0)=P_Nx,\label{e10}\\[1mm]
\dd\frac{dX_s}{dt}+\cala_sX_s=0,\quad X_s(0)=(I-P_N)x.\label{e11}
\end{eqnarray}
We set $A_u=\left\{\<\cala\vf_j,\vf^*_k\>\right\}^N_{j,k=1}={\rm
diag}\|\lbb_j\|^N_{j=1}$ and so, by \eqref{e8}, we may rewrite
\eqref{e10} in terms of
$y=\left\{y_j=\<X_u,\vf^*_j\>\right\}^N_{j=1}$ as
\begin{equation}\label{e12}
\frac{dy}{dt}+A_uy=0,\quad
y(0)=\left\{\<P_Nx,\vf^*_j\>\right\}^N_{j=1}.\end{equation} Since
${\rm Tr}(-A_u)=-\lbb_1-\lbb_2-\cdots-\lbb_N<0$, it follows by
Theo\-rem~2 in \cite{3} that there is a sequence of skew-symmetric
matrices $\{C^k\}^M_{k=1}$, where $M=N-1$ such that the solution
$y$ to the Stratonovich stochastic system
\begin{equation}\label{e13}
dy+A_uy\,dt=\sum^M_{k=1}C^ky\circ d\b_k,\quad t\ge0,\end{equation}
has the property
\begin{equation}\label{e14}
|y(t)|\le C|y(0)|e^{-\g_0 t},\ \pas, \ \ff t>0,\end{equation}
where $\g_0>0.$  The matrix $C^k$ is explicitly constructed in
\cite{3} and it will be used below to construct a stabilizable
feedback controller of the form \eqref{e3}. Namely, we set in
\eqref{e3}
\begin{equation}\label{e15}
R_k(X)=\sum^N_{i,j=1}C^k_{ij}\<X,\vf^*_j\>\phi_i,\quad k=1,...,M,\end{equation}
where $\|C^k_{ij}\|^N_{i,j=1}=C^k$,
\begin{equation}\label{e16}
\phi_i=\sum^N_{\ell=1}\a_{i\ell}\vf^*_\ell,\quad i=1,...,N,\end{equation}
and $\a_{i\ell}$ are chosen in such a way that
\begin{equation}\label{e17}
\sum^N_{\ell=1}\a_{i\ell}\g_{\ell j}=\de_{ij},\quad
i,j=1,...,N.\end{equation} Here, $\g_{\ell
j}=\int_{\calo_0}\vf^*_\ell\ov\vf^*_jd\xi$ and  since, by the
unique continuation property (see \cite{4}, p.~157), the
eigenfunction system $\{\vf^*_j\}$ is li\-nearly independent on
$\calo_0$,  we infer that the matrix $\|\g_{\ell
j}\|^N_{\ell,j=1}$ is not singular and, therefore, there is a
unique system $\{\a_{i\ell}\}$ which satisfies \eqref{e17}. Then,
by \eqref{e16}, we see that
\begin{equation}\label{e18}
\<{\bf1}_{\calo_0}\phi_i,\vf^*_j\>=\de_{ij},\quad
i,j=1,...,N.\end{equation} As mentioned earlier, $\calo_0$ is an
arbitrary open subset of $\calo$.

Theorem \ref{t1} is the main result.

\begin{theorem}\label{t1} The solution $X$ to the closed loop system \eqref{e4}, where $R_k$ are defined by \eqref{e15}, is exponentially stable with probability one, that is,
\begin{equation}\label{e19}
P\left[|X(t)|_{\calh}\le Ce^{-\g t}|x|_H,\quad\ff t\ge0\right]=1,\end{equation}
where $\g>0$ and $C>0$ are independent of $\omega\in\Omega$.
\end{theorem}

\subsection{The proof of Theorem \ref{t1}}

The idea of the proof, already used in stabilization theory of
infinite-dimensional systems with finite-dimensional unstable
subspaces, is to stabilize the finite-dimensional system
\eqref{e10} by a feedback controller of the form \eqref{e3} and to
reconstruct consequently the system via the infinite-dimensional
stable complement \eqref{e11}.

Namely, taking into account \eqref{e10}, \eqref{e11}, we write the
closed loop system \eqref{e4}, that is,
\begin{equation}\label{e20}
dX+\cala X\,dt=P\left[{\bf1}_{\calo_0}\sum^M_{k=1}R_k(X)\circ
d\b_k\right]\end{equation} as
\begin{equation}\label{e21}
dX_u+\cala_uX_udt=P_N\sum^M_{k=1}\sum^N_{i,j=1}C^k_{ij}\<X_u,\vf^*_j\>
P(\phi_i)\circ d\b_k,\end{equation}
\begin{equation}\label{e22}
dX_s+\cala_sX_sdt=(I-P_N)\sum^M_{k=1}
\sum^N_{i,j=1}C^k_{ij}\<X_u,\vf^*_j\>P(\phi_i)\circ d\b_k,
\end{equation}
where $X_u+X_s=X$.

Taking into account \eqref{e18} and that
$X_u=\sum^N_{j=1}y_j\vf_j$, we may rewrite \eqref{e21} as
\begin{equation}\label{e23}
dy_\ell+\lbb_\ell y_\ell dt=\sum^M_{k=1}C^k_{\ell j}y_j\circ d\b_k,\quad\ell=1,...,N,\end{equation}
and so, by \eqref{e14}, we have that
\begin{equation}\label{e24}
|y_\ell(t)|\le Ce^{-\g t}|y_\ell(0)|,\ \pas,\ \ff t>0,\ \ell=1,...,N.\end{equation}
(We have denoted by $C$ several positive constants independent of $t$ and $\omega\in\Omega$.)

As regards the existence of a solution $X_s$ to \eqref{e22}, this
is standard and follows from the general theory of linear
infinite-dimensional stochastic equations.

 Now, in order to
estimate $X_s$, it is convenient to replace \eqref{e22} by its Ito
formulation (see, e.g., \cite{9})
\begin{equation}\label{e25}
\barr{rcl} dX_s+\cala_sX_s\,dt&=&\dd\frac12\
\sum^M_{k=1}[P({\bf1}_{\calo_0}R_k]^2X_udt\vsp
&+&\dd\sum^M_{k=1}P({\bf1}_{\calo_0}R_k(X_u))d\b_k.\earr\end{equation}
Taking into account that,  by \eqref{e7}, $e^{-\cala_st}$ is an
exponentially stable semigroup on $\calx_s$, without loss of
generality we may assume that $\Res\,\<\cala_sx,x\>\ge\g|x|^2_H$.
(Otherwise, proceeding as in \cite{5}, we replace the scalar
product $\<x,y\>$ by $\<Qx,y\>$ where $Q$ is the solution to the
Lyapunov equation $\cala_sQ+Q\cala^*_s=\g I$.) Then, applying
Ito's formula in \eqref{e25}, we see that
\begin{equation}\label{e26}
\barr{l} d|X_s(t)|^2_\calh{+}2\g|X_s(t)|^2_\calh
dt{=}\<X_s(t),\dd\sum^M_{k=1}[P({\bf1}_{\calo_0}R_k]^2X_u(t)\>dt\vsp
\qquad+\dd\sum^M_{k=1} |P({\bf1}_{\calo}R_k(X_u))|^2dt
+2\dd\sum^M_{k=1}P({\bf1}_{\calo}R_k(X_u))X_sd\b_k.
\earr\hspace*{-3mm}\end{equation} Taking into account that, by
\eqref{e24}, $|X_u(t)|_\calh\le e^{-\g t}|X_u(0)|,$ we infer by
\eqref{e26} that
$$E|X_s(t)|^2_\calh\le C|X_s(0)|^2e^{-\g t},\quad\ff t\ge0.$$
Since
$M(t)=\int^t_0\sum^M_{k=1}P({\bf1}_{\calo_0}R_k(X_u)X_s\,d\b_k$ is
a local martingale, it follows by \eqref{e26} and Lemma 3.1 in
\cite{5} that
$$|X_s(t)|_\calh\le C|X_s(0)|e^{-\g t},\quad\ff t\ge0,\ \pas,$$
which completes the proof.

\section{The design of a real valued noise controller}

A nice feature of the feedback controller given by
Theorem~\ref{t1} is its simple structure. Moreover, its
computation does not rise any special problem because the
stabilizing matrices $C^k$ for the diagonal system \eqref{e23} can
be explicitly expressed.

If all $\lbb_j$, $j=1,...,N$, are real, then the feedback
controller $\{R_k\}$ given by Theorem \ref{t1} as well as the
closed loop system \eqref{e4} are real, too. However, if some
$\lbb_j$ are complex, then the above feedback controller does not
stabilize system \eqref{e1}, but its complex realization in the
space $\calh=H+iH$. In this case, though the stabilizing feedback
controller \eqref{e3} is quite simple, its implementation into
real systems rises some delicate problems because  the feedback
controller is in implicit form as function of $\Res\,X$,
$\Ims\,X$. In order to circumvent these inconvenience, we will
derive below from the above construction a real valued noise
feedback controller which has a stabilizing effect on the
Oseen--Stokes system \eqref{e1}.

To this end, we denote by $\{\psi_j\}^N_{j=1}$ the  orthogonalized
system (via Schmidt procedure) $\{\Res\,\vf_j,\Ims\,\vf_j,\
j=1,...,N\}$. Taking into account that each $\vf_j$ in the above
system arises together with its conjugate, the dimension of the
system $\{\psi_j\}$ is again $N$. We also set
$$\calx_u^{\rm re}={\rm lin\ span}\{\psi_j\}^N_{j=1}$$
and note that  $\calx_u^{\rm re}=\{\Res\,y;\ y\in\calx_u\}$ and
$H=\calx^{\rm re}_u\oplus\calx_s^{\rm re}$, where $X^{\rm
re}_s=\{\Res\,y;\ y\in\calx_s\}$.    Moreover, the operator
$\cala_0$ leaves in\-va\-riant spaces $\calx^{\rm re}_u,$
$\calx_s^{\rm re}$ and if we set
$$\cala^{\rm re}_u=\cala_0|_{\calx_u^{\rm re}},
\quad\cala_s^{\rm re}=\cala_0|_{\calx^{\rm re}_s},$$ we have, of
course,
\begin{equation}\label{e27}
\|e^{-\cala^{\rm re}_s t}\|_{L(H)}\le Ce^{-\g_0 t},\quad\ff t>0,\end{equation}
for some $\g_0>0$. It is also easily seen that $${\rm Tr}[-\cala^{\rm re}_u]={\rm Tr}[-\cala_u]<0.$$
Then, we define as in \eqref{e15} the operators
\begin{equation}\label{e28}
\wt R_k(x)=\sum^N_{i,j=1}\wt C^k_{ij}\<X,\vf_j\>\wt\phi_i,\quad
k=1,....,M,\end{equation} where $\wt C^k=\|\wt
C^k_{ij}\|^N_{i,j=1}$  is the matrix system which stabilizes with 
probability one, via Theorem 2 in \cite{3}, the finite-dimensional
system
$$dy+A^{\rm re}_uy\,dt=\sum^M_{k=1}\wt C^k y\circ d\b_k.$$
Here, $A^{\rm re}_u=\|\<\cala^{\rm
re}_u\psi_i,\psi_j\>\|^N_{i,j=1}$  and $\wt\phi_i$ is given by
$$\wt\phi_i=\sum^N_{\ell=1}\wt\a_{i\ell}\psi_\ell,\quad i=1,....,N,$$
where $\wt\a_{i\ell}$ are chosen in such a way that
$$\sum^N_{\ell=1}\wt\a_{i\ell}\wt\g_{\ell j}=\de_{ij},\quad\wt\g_{\ell j}=\int_{\calo_0}\psi_\ell\psi_j d\xi.$$
As in the previous case, the system $\{\psi_j\}^N_{j=1}$  is still
linearly independent on $\calo_0$ which implies that
$\det\|\wt\g_{\ell j}\|\ne0$ and the above system has a unique
solution $\wt\a_{i\ell}$.

Then, arrived at this point,  the proof of Theorem \ref{t1}
applies neatly to conclude that the feedback noise controller
\eqref{e3} with $\wt R_k$ instead of $R_k$ is stabilizing in the
probability system \eqref{e9}.

Namely, we have the following stabilization result.

\begin{theorem}\label{t2} The solution $X$ to the closed loop
system \eqref{e4} with $R_k=\wt R_k$ defined by \eqref{e28}
satisfies \eqref{e19} in the real norm $|X(t)|_H$.\end{theorem}

\begin{remark} {\rm It should be emphasized that there is a close
connection  between the unique continuation property of
eigenfunctions $\vf_j$ (or $\vf^*_j$) of the Oseen--Stokes
operator and the above construction of a stabilizing
finite-dimensional feedback controller \eqref{e3}. In fact, as
seen above, the design of $u$ in the form \eqref{e15} or
\eqref{e28} is essentially based on this sharp
property through existence in the algebraic system \eqref{e17}.}\end{remark}

It should be mentioned also that, likewise the cor\-res\-pon\-ding Ito noise con\-trol\-ler (see \cite{4}), the stabilizable feedback con\-trol\-ler \eqref{e3} is robust.

\begin{remark} {\rm One might speculate that the noise feedback controller
\begin{equation}\label{e29}
u={\bf1}_{\calo_0}\ \sum^M_{k=1}R_k(X-X_e)\circ
d\b_k\end{equation} inserted in the right hand side of
Navier--Stokes system \eqref{e2} stabilizes exponentially in
probability the equilibrium solution $X_e$ to \eqref{e2}. In general, this might not be true, but it
happens for Ito noise of the form \eqref{e29} for any $x$ in a sufficiently small neighborhood of $X_e$ (see \cite{7}) and
one might expect that the fixed point argument used there for the
equivalent random system is still applicable in the present case.
We expect to give details in a later work.}\end{remark}

\section{Conclusions}

We have designed in this paper a  Stratonovich stochastic feedback
controller which exponentially stabilizes in pro\-ba\-bi\-lity a
general Oseen--Stokes system from fluid dynamics. The controller
has the support in an arbitrary open subset $\calo_0$ of the
velocity field  domain $\calo\subset\mathbb{R}^d$, $d=2,3,$ and
has a finite-dimensional linear structure which involves the dual
eigenfunctions corresponding to unstable eigenvalues of the
systems. The stabilization effect is independent of the Reynold
number $1/\nu$ through the dimension $N$ of the stabilizing controller \eqref{e3} might depend on $\nu$.

\bigskip\n{\bf Acknowledgement.} This work is
supported by CNCSIS project PN II IDEI ID\_70/2008.

\end{document}